\documentclass[12pt, leqno]{amsart}

\usepackage[OT2,T1]{fontenc}
\usepackage{indentfirst}
\usepackage{amstext}
\usepackage{amsthm}
\usepackage{amsopn}
\usepackage{amsfonts}
\usepackage{amsmath}
\usepackage{latexsym}
\usepackage[francais,english]{babel}
\usepackage{amscd}
\usepackage{amssymb}
\usepackage{amsmath}
\usepackage[all,cmtip]{xy}
\usepackage{leftidx}
\usepackage{graphicx}
\usepackage{etoolbox}
\patchcmd{\section}{\normalfont\scshape\centering}{\normalfont\bfseries}{}{}
\patchcmd{\subsection}{-.5em}{.5em}{}{}

\newtheorem{theo}{{Theorem}}[section]
\newtheorem{coro}[theo]{{Corollary}}

\newtheorem{prop}[theo]{Proposition}

\theoremstyle{definition}

\newtheorem{defn}[theo]{Definition}
\newtheorem{example}[theo]{Example}
\numberwithin{equation}{section}

\newtheorem{notation}[theo]{Notation}

\DeclareFontEncoding{OT2}{}{} 
  \newcommand{\textcyr}[1]{%
    {\fontencoding{OT2}\fontfamily{wncyr}\fontseries{m}\fontshape{n}%
     \selectfont #1}}
\newcommand{\sha}{{\mbox{\textcyr{Sh}}}}

\begin{document}
\tolerance 400 \pretolerance 200 \selectlanguage{english}

\title{Alexander polynomials and signatures of some high-dimensional knots}
\author{Eva  Bayer-Fluckiger}
\date{\today}
\maketitle

\begin{abstract}
We give necessary and sufficient conditions for an integer to be the signature of a $4q-1$-knot in $S^{4q+1}$ 
with a given square-free Alexander polynomial. 

\medskip

\end{abstract}

\small{} \normalsize

\medskip

\selectlanguage{english}
\section{Introduction} What are the possiblities for the signatures of knots with a given Alexander polynomial ? This question is answered in \cite{B 21} for ``classical" knots, i.e. knots $K^1 \subset S^3$, with some
restrictions on the Alexander polynomial, and the same
results hold for knots $K^m \subset S^{m+2}$ if $m  \equiv 1 \ {\rm (mod \ 4)}$. In the present paper, we
consider high-dimensional knots $K^m \subset S^{m+2}$ with $m  \equiv -1 \ {\rm (mod \ 4)}$. In the introduction,
we describe the results for $m > 3$; the case
$m = 3$ is somewhat different (see Section \ref{3}).

\medskip
Let $m \geqslant 7$ be an integer with  $m  \equiv -1 \ {\rm (mod \ 4)}$. An {\it $m$-knot} $K^m \subset S^{m+2}$
is by definition a smooth, oriented submanifold of $S^{m+2}$, homeomorphic to $S^m$; in the following,
a {\it knot} will mean an $m$-knot as above. We refer to the book of Michel and Weber \cite {MW} for a survey
of high-dimensional knot theory. 


\medskip 
Let $K^m$ be a knot. The {\it Alexander polynomial} $\Delta = \Delta_K \in {\bf Z}[X]$
is a polynomial of even degree; set $2n = {\rm deg}(\Delta)$. It satisfies  the following three properties (see for instance \cite {Le 69}, Proposition 1) :

\medskip
(1) $\Delta(X) = X^{2n} \Delta(X^{-1})$,

\medskip
(2) $\Delta(1) = (-1)^n$,

\medskip
(3) $\Delta(-1)$ is a square.

\medskip
Conversely, if $\Delta \in {\bf Z}[X]$ is a degree $2n$ polynomial satisfying conditions (1)-(3), then there
exists a knot with Alexander polynomial $\Delta$ (cf. Levine \cite {Le 69}, Proposition 2 and Lemma 3; note
that Lemma 3 is based on a result of Kervaire, \cite {K 65}, Th\'eor\`eme II.3).

\medskip
Let $F^{m+1}$ be a Seifert
hypersurface of $K^m$ (see for instance \cite{MW}, Definition 6.16), and
let $L = H_n(F^{m+1},{\bf Z})/{\rm Tors}$, where ${\rm Tors}$ is the ${\bf Z}$-torsion subgroup of 
$H_n(F^{m+1},{\bf Z})$. Let $S : L \times L \to {\bf Z}$ be the intersection form; since $m  \equiv -1 \ {\rm (mod \ 4)}$, the form $S$ is {\it symmetric}. 
The {\it signature} of $K^m$ is by definition the signature of the symmetric form $S$; it is an invariant
of the knot. The form $S$ is
even and unimodular, therefore its signature is $ \equiv 0 \ {\rm (mod \ 8)}$.

\medskip Levine's construction (see  \cite {Le 69}, Proposition 2 and Lemma 3) shows the existence
of a knot with Alexander polynomial $\Delta$ and signature $0$. It is natural to ask : {\it what other
signatures occur}~?

\medskip Let us denote by $\rho (\Delta)$ the number of roots $z$ of $\Delta$ such that $|z| = 1$. If a
knot has Alexander polynomial $\Delta$ and signature $s$, then $|s| \leqslant \rho (\Delta)$. 
This shows that the conditions $s \equiv 0 \ {\rm (mod \ 8)}$ and $|s| \leqslant \rho (\Delta)$
are {\it necessary} for the existence of a knot with Alexander polynomial $\Delta$ and signature $s$; however,
these conditions are not sufficient, as shown by the following example, taken from \cite{GM}, Proposition 5.2 :

\medskip
\noindent
{\bf Example 1.} Let $\Delta(X) = (X^6 - 3 X^5 - X^4 + 5 X^3 - X^2 - 3X + 1)(X^4 - X^2 + 1)$; we have
$\rho(\Delta) = 8$, hence $s = -8,0$ and $8$ satisfy the above necessary conditions. However, there does
not exist any knot with Alexander polynomial $\Delta$ and signature $-8$ or $8$.

\medskip
Let $\Delta \in {\bf Z}[X]$ be a polynomial satisfying conditions (1)-(3), and suppose that $\Delta$ is {\it square-free}. We associate to $\Delta$ a finite abelian group $G_{\Delta}$ that controls the signatures of the knots
with Alexander polynomial $\Delta$ (see \S \ref{group} - \S  \ref{knot section}). In particular, we have (cf. Corollary \ref{knot coro sign}) : 

\medskip
\noindent
{\bf Theorem 1.} Assume that  $G_{\Delta} = 0$, and let $s$ be an integer with $s \equiv 0 \ {\rm (mod \ 8)}$ and $|s| \leqslant \rho(\Delta)$. Then there exists a knot with Alexander polynomial $\Delta$ and signature $s$.

\medskip
The vanishing of the group $G_{\Delta}$ has other geometric consequences : 
we show the existence of {\it indecomposable} knots with Alexander polynomial $\Delta$ (see \S \ref{indecomposable}).

\section{Seifert forms and Seifert pairs}\label{Seifert}

Seifert forms  are well-known objects of knot theory; the aim of this section is
to recall this notion, and to show that it is equivalent to the one of {\it Seifert pairs}; this notion was
introduced, under a different name, by Kervaire in \cite {K 71} in the context of knot cobordism; 
see also Stoltzfus (\cite {St}) and \cite{B 82}, \S 5.

\begin{defn} \label{Seifert form}
A {\it Seifert form} is by definition a pair $(L,A)$, where $L$ is a free $\bf Z$-module of finite rank
and $A : L \times L \to {\bf Z}$ is a $\bf Z$-bilinear form such that the symmetric form $L \times L \to {\bf Z}$
sending $(x,y)$ to $A(x,y) + A(y,x)$ is unimodular (i.e. has determinant $\pm 1$); the  {\it signature}
of $(L,A)$ is by definition the signature of this symmetric form. 

\medskip The {\it Alexander polynomial} of $(L,A)$, denoted by $\Delta_A$,
is by definition the determinant of the form $L \times L \to {\bf Z}[X]$ given by
$$(x,y) \mapsto A(x,y)X + A(y,x).$$
\end{defn}

\begin{defn} \label{Seifert pair}
A {\it Seifert pair} is by definition a triple $(L,S,a)$, where $L$ is a free $\bf Z$-module of finite rank,
$S : L \times L \to {\bf Z}$ is an even (i.e. $S(x,x)$ is an even integer for all $x \in L$), unimodular, symmetric $\bf Z$-bilinear form, and $a : L \to L$ is an injective $\bf Z$-linear map such that
$$S(ax,y) = S(x,(1-a)y)$$ for all $x,y \in L$.
\end{defn}

\medskip
Let $(L,S,a)$ be a Seifert pair. Since $S$ is even and unimodular, the rank of $L$ is an even integer; let
$n \in {\bf Z}$ be such that ${\rm rank}(L) = 2n$. Let $A : L \times L \to {\bf Z}$ be defined by
$$A(x,y) = S(ax,y);$$  note that $(L,A$) is a Seifert form, and we have

\begin{prop}\label{bijection} Sending $(L,S,a)$ to $(L,A)$ as above induces a bijection between
isomorphism classes of Seifert pairs and of Seifert forms. Let $P_a$ be the characteristic polynomial
of $a$. We have $$P_a(X) = (-1)^n X^{2n} \Delta_A(1 - X^{-1}).$$
\end{prop}

\medskip

Note that $\Delta_A(X) = X^{2n} \Delta_A(X^{-1})$, and that $P_a(X) = P_a(1-X)$. 

\begin{defn}
A {\it lattice} is a pair $(L,S)$, where $L$ is a free $\bf Z$-module of finite rank, and $S : L \times L \to {\bf Z}$
is a symmetric bilinear form with ${\rm det}(S) \not = 0$. We say that $(L,S)$ is unimodular if ${\rm det}(S) = \pm 1$, and even if $S(x,x)$ is an even integer for all $x \in L$.
\end{defn}

Note that a Seifert pair consists of an even, unimodular lattice $(L,S)$ and an injective endomorphism
$a : L \to L$ such that $S(ax,y) = S(x,(1-a)y)$ for all $x,y \in L$,

\section {Involutions of $K[X]$, symmetric polynomials and bilinear forms compatible with a module}\label{Witt}

Let $K$ be a field, let $R$ be a commutative $K$-algebra, and let $\sigma : R \to R$ be an involution; we say that $\lambda \in R$ is
{\it $\sigma$-symmetric} (or symmetric, if the choice of $\sigma$ is clear from the context) if $\sigma(\lambda) = 
\lambda$.

\begin{example} (1)  Let $\sigma : K[X,X^{-1}] \to K[X,X^{-1}]$ be the involution sending $X$ to $X^{-1}$. If $\Delta$  is the Alexander polynomial of a Seifert form of rank $2n$, then $X^{-n} \Delta(X)$
is symmetric.

\medskip
(2) Let $\sigma : K[X] \to K[X]$ be the involution sending $X$ to $1-X$; the symmetric polynomials are the
$f \in K[X]$ such that $f(1-X) = f(X)$. The characteristic polynomial of a Seifert pair is symmetric.

\end{example} 

Let $M$ be an $R$-module that is a finite dimensional $K$-vector space. Recall from \cite {B 21}, \S 1, that a non-degenerate symmetric bilinear form
$b : M \times M \to K$ is called an $(R,\sigma)$-{\it bilinear form} if $$b(\lambda x,y) = b(x, \sigma(\lambda) y)$$
for all $x,y \in V$ and for all $\lambda \in R$. 

\begin{example} Let $(L,S,a)$ be a Seifert pair and set $V = L \otimes_{\bf Z}{\bf Q}$; we denote
by $S : V \times V \to {\bf Q}$ and  $a : V \to V$ the symmetric bilinear form and the $\bf Q$-linear
map induced by $S$ and $a$. Let $\sigma :
{\bf Q}[X] \to {\bf Q}[X]$ be the involution sending $X$ to $1-X$. We endow $V$ with a structure
of ${\bf Q}[X]$-module by setting $X.x = a(x)$ for all $x \in V$; note that $S : V \times V \to {\bf Q}$
is a $({\bf Q}[X],\sigma)$ bilinear form. 

\end{example}

Let $V$ be a finite dimensional $K$-vector space, and let $q : V \times V \to K$ be a non-degenerate
symmetric bilinear form. Following \cite {B 21}, \S 1, we say that $M$ and $(V,q)$ are {\it compatible}
if there exists a $K$-linear isomorphism $\phi : M \to V$ such that the bilinear form
$b_{\phi} : M \times M \to K$, defined by $b_{\phi}(x,y) = q(\phi(x),\phi(y))$, is an $R$-bilinear form. 

\begin{example}  Let $\sigma : K[X] \to K[X]$ be the involution sending $X$ to $1-X$, and let $P \in K[X]$
be a monic,  $\sigma$-symmetric polynomial. Assume $P$ is a product of distinct monic, symmetric,
irreducible factors; let us denote by $I$ the set of these polynomials. Set $M = \underset{f \in I} \oplus
K[X]/(f)$, regarded as a $K[X]$-module.
Let $V$ be a finite dimensional $K$-vector space, and let $q : V \times V \to K$ be a non-degenerate
symmetric bilinear form. The module $M$ and the form $(V,q)$ are compatible if and only if there
exists an endomorphism $a : V \to V$ with characteristic polynomial $P$ such that $q(ax,y) = q(x,(1-a)y)$
for all $x,y \in V$. 

\end{example}

\section{Milnor signatures}\label{Milnor}

We recall the notion of {\it Milnor signatures}, introduced by Milnor  in \cite {M 68}, in the context of Seifert pairs.
Let $(L,S,a)$ be a Seifert pair, and let $P \in {\bf Z}[X]$ be the characteristic polynomial of $a$. Assume
that the polynomial $P$ is {\it square-free}, i.e. has no repeated factors; we also suppose that if
$f \in {\bf Z}[X]$ is a monic, irreducible factor of $P$, then $f(X) = f(1-X)$. 

\medskip Let $V = L \otimes_{\bf Z}{\bf R}$. Let $f \in {\bf R}[X]$ be a monic, irreducible factor
of degree $2$ of $P \in {\bf R}[X]$; note that this implies that $f(X) = f(1-X)$.

\begin{defn} The {\it signature} of $(L,S,a)$ at $f$ is by definition the signature of the restriction of
$S$ to ${\rm Ker}(f(a))$. 

\end{defn}

\begin{notation} Let ${\rm Irr}_{\bf R}(P)$ be the set of monic, irreducible factors $f \in {\bf R}[X]$  of
degree $2$ of $P$. Let $s \in {\bf Z}$. We denote by ${\rm Mil}(P)$ the set of maps
$${\rm Irr}_{\bf R}(P) \to \{-2,2 \},$$ and by ${\rm Mil}_s(P)$ the set of $\tau \in {\rm Mil}(P)$ such that

$$\underset {f \in {\rm Irr}_{\bf R}(P)} \sum \tau(f) = s.$$

\end{notation}

Let $n \geqslant 1$ be an integer, and let $\Delta \in {\bf Z}[X]$ be a polynomial of degree $2n$ such
that $\Delta(X) = X^{2n}\Delta(X^{-1})$, $\Delta(1) = (-1)^n$ and that $\Delta(-1)$ is a square of an
integer. Suppose that  $P(X) = (-1)^nX^{2n} \Delta(1-X^{-1})$. We define ${\rm Mil}_s(\Delta)$
as in \cite {B 21}, \S 26; note that there are obvious bijections between ${\rm Irr}_{\bf R}(P)$ and
${\rm Irr}_{\bf R}(\Delta)$, ${\rm Mil}_s(P)$ and ${\rm Mil}_s(\Delta)$, and that we recover the
usual notion of Milnor signature.

\medskip
If $P$ and $\Delta$ are as above, set $\rho(P) = \rho(\Delta)$; alternatively, $\rho(P)$ can be
defined as the number of roots $z$ of $P$ with $z + \overline z = 1$, where $\overline z$ denotes
the complex conjugate of $z$. Note that $\rho(\Delta) = |{\rm Irr}_{\bf R}(\Delta)|$ and 
$\rho(P) = |{\rm Irr}_{\bf R}(P)|$.

\section{The obstruction group}\label{group}

Let $P \in {\bf Z}[X]$ be a monic polynomial such that $P(1-X) = P(X)$. Assume that $P$ is a product
of distinct irreducible monic polynomials $f \in {\bf Z}[X]$ such that $f(1-X) = f(X)$. We associate to
$P$ an elementary abelian $2$-group $G_P$ that will be useful in the following sections; this construction
is similar to the one of \cite{B 21}, \S 21. 

\medskip Let $I$ be the set of irreducible factors of $P$. If $f,g \in I$, let $\Pi_{f,g}$ be the set of
prime numbers $p$ such that $f \ {\rm mod} \ p$ and $g \ {\rm mod} \ p$ have a common factor
$h \in {\bf F}_p[X]$ such that $h(1-X) = h(X)$. Let $C(I)$ be the set of maps $I \to {\bf Z}/2{\bf Z}$, and
let $C_0(I)$ be the set of $c \in C_0(I)$ such that $c(f) = c(g)$ if  $\Pi_{f,g} \not = \varnothing$. Note
that $C_0(I)$ is a group with respect to the addition of maps, and let $G_P$ be the quotient of
the group $C_0(I)$ by the subgroup of the constant maps. 

\begin{example}\label{first example} Let $f_1(X) = X^4 - 2X^3 + 5X^2 - 4X + 1$ and $$f_2(X) = X^4 - 2X^3 + 11X^2 - 10X + 3;$$ set
$P = f_1 f_2$. We have $\Pi_{f_1,f_2} = \{2\}$, hence $G_P = 0$. 

\end{example}

\medskip If $P(X) = (-1)^nX^{2n} \Delta(1-X^{-1})$ for some polynomial $\Delta \in {\bf Z}[X]$, 
set $G_{\Delta} = G_P$. If moreover $\Delta(0) = \pm 1$, then the
group $G_{\Delta}$ is equal to the obstruction group $\sha_{\Delta(0) \Delta}$ of \cite{B 21}, \S 21 and \S 25. In particular,
25.8 - 25.11, 31.4 and 31.5  of \cite {B 21} provide examples of obstruction groups in our context as well. This is
also the case for the following example, given in the introduction~:

\begin{example} Let $g_1(X) = X^6 - 3 X^5 - X^4 + 5 X^3 - X^2 -3X +1$ and 
$g_2(X) = X^4 - X^2 + 1$; set $\Delta = g_1 g_2$, as in Example 1. Set $f_1(X) = -X^6g_1(1-X^{-1})$,
and $f_2 (X) = X^4 g_2(1-X^{-1})$, and let $P = f_1 f_2$. The polynomials $f_1$ and $f_2$ are relatively
prime over $\bf Z$, hence $\Pi_{f_1,f_2} = \varnothing$; therefore $G_P \simeq {\bf Z}/2{\bf Z}$. 

\end{example}

\section{Seifert pairs with a given characteristic polynomial and signature}\label{obstruction}

Let $n \geqslant 1$ be an integer, and let $\Delta \in {\bf Z}[X]$ be a polynomial of degree $2n$ such
that $\Delta(X) = X^{2n}\Delta(X^{-1})$, $\Delta(1) = (-1)^n$ and that $\Delta(-1)$ is a square of an
integer. Set $P(X) = (-1)^nX^{2n} \Delta(1-X^{-1})$. Assume that $P$ is a product
of distinct irreducible monic polynomials $f \in {\bf Z}[X]$ such that $f(1-X) = f(X)$, and let
$I$ be the set of irreducible, monic factors of $P$. 

\medskip Let $G_P$ be the group introduced in \S \ref{group}, and set $G_{\Delta} = G_P$. 

\medskip Let $s$ be an integer such that $s \equiv 0 \ {\rm (mod \ 8)}$, and that $|s| \leqslant \rho(P)$. 
Let $\tau \in {\rm Mil}_s(P)$. The aim of this section is to give a necessary and sufficient condition
for the existence of a Seifert pair with characteristic polynomial $P$ and Milnor signature $\tau$. 

\medskip
Let $V$ be a $\bf Q$-vector space of dimension $2n$, and let $S : V \times V \to {\bf Q}$ be a
non-degenerate quadratic form of signature $s$ containing an even, unimodular lattice; such a
form exists and is unique up to isomorphism (see for instance \cite {B 21}, Lemma 25.5). 

\medskip
Let $M = \underset{f \in I} \oplus {\bf Q}[X]/(f)$, considered as a ${\bf Q}[X]$-module. Let
$\sigma : {\bf Q}[X] \to {\bf Q}[X]$ be the $\bf Q$-linear involution such that $\sigma(X) = 1-X$. 
The Milnor signature $\tau \in {\rm Mil}_s(P)$ determines an $({\bf R}[X],\sigma)$-quadratic form (cf \cite {B 21},
Example 24.1). The local conditions of \cite {B 21}, \S 24 are satisfied. Indeed, the ${\bf R}[X]$-module
$M \otimes {\bf R}$ is compatible with $(V,S)$ by \cite {B 15}, Proposition 8.1. Using a result
of Levine (see \cite{Le 69}, Proposition 2) and the bijection between Seifert forms and Seifert pairs
(see \S \ref{Seifert}), we see that there exists a Seifert pair of characteristic polynomial $P$. This implies
that for all prime numbers $p$, the ${\bf Q}_p[X]$-module $M \otimes_{\bf Q} {\bf Q}_p$ and the quadratic form
$(V,S) \otimes _{\bf Q} {\bf Q}_p$ are compatible. 

\medskip
As in \cite{B 21}, \S 24, we define a homomorphism $\epsilon_{\tau} : G_P \to {\bf Z}/2{\bf Z}$. 

\begin{theo}\label{main} There exists a Seifert pair with characteristic polynomial $P$ and Milnor signature $\tau$
if and only if $\epsilon_{\tau} = 0$.

\end{theo}

\noindent
{\bf Proof.} By \cite {B 21}, Theorem 24.2, the global conditions are satisfied if and only if $\epsilon_{\tau} = 0$. 
Using \cite {B 21}, Proposition 6.2 this is equivalent with the existence of a Seifert pair having characteristic
polynomial $P$ and Milnor signature $\tau$. 

\begin{coro}\label{main coro} Assume that $G_P = 0$. Then for all $\tau \in {\rm Mil}_s(P)$ there
exists a Seifert pair with characteristic polynomial $P$ and Milnor signature $\tau$.

\end{coro} 

\section{Seifert forms with a given Alexander polynomial and signature}\label{Seifert form section}

We keep the notation of the previous section. Using Proposition \ref{bijection}, Theorem \ref{main} and Corollary \ref{main coro} can be reformulated as follows :

\begin{theo}\label{main Seifert}  There exists a Seifert form with Alexander polynomial $\Delta$ and Milnor signature $\tau$
if and only if $\epsilon_{\tau} = 0$.

\end{theo}

\begin{coro}\label{main coro Seifert} Assume that $G_{\Delta} = 0$. Then for all $\tau \in {\rm Mil}_s(\Delta)$ there
exists a Seifert form with Alexander polynomial $\Delta$ and Milnor signature $\tau$.

\end{coro}

\section{Knots with a given Alexander polynomial and signature}\label{knot section}

We keep the notation of the previous two sections. Let $m \geqslant 7$ be an integer with  $m  \equiv -1 \ {\rm (mod \ 4)}$. We refer to \cite {MW}, 6.5 for the definition of the Seifert form associated to an $m$-knot. The results of this section rely on a result of Kervaire~:

\begin{theo} \label{Kervaire} Let $(L,A)$ be a Seifert form. Then there exists an $m$-knot with associated Seifert
form isomorphic to $(L,A)$. 

\end{theo}

\noindent 
{\bf Proof.} This is proved by Kervaire in \cite {K 65}, Theorem II.3, and formulated more explicitly
by Levine in \cite{Le 69}, Lemma 3 and \cite{Le 70}, Theorem 2. A different proof is given by Michel
and Weber in \cite {MW}, Theorem 7.3 (see also the remark at the end of \cite {MW}, \S 7.1). 

\medskip

Combining Theorem \ref{main Seifert} and Corollary \ref{main coro Seifert} with Theorem \ref{Kervaire},
we 
have the following applications :

\begin{theo}\label{knot} There exists an $m$-knot with Alexander polynomial $\Delta$
and Milnor signature $\tau$ if and only if $\epsilon_{\tau} = 0$. 

\end{theo} 

\begin{coro}\label{knot coro} Assume that $G_{\Delta} = 0$. Then for all $\tau \in {\rm Mil}_s(\Delta)$ there
exists an $m$-knot with Alexander polynomial $\Delta$ and Milnor signature $\tau$.

\end{coro}

Recall that $s$ is an integer such that $s \equiv 0 \ {\rm (mod \ 8)}$, and that $|s| \leqslant \rho(\Delta)$. 

\begin{coro}\label{knot coro sign}  Assume that $G_{\Delta} = 0$. Then there
exists an $m$-knot with Alexander polynomial $\Delta$ and signature $s$.

\end{coro}

\section {Indecomposable knots with decomposable Alexander polynomial}\label{indecomposable}

As an application of Corollary \ref{knot coro sign}, we give some examples of indecomposable
knots with decomposable Alexander polynomials. Let $m \geqslant 7$ be an integer with  $m  \equiv -1 \ {\rm (mod \ 4)}$. 

\begin{example}\label{4} Let $\Delta = \Delta_1 \Delta_2$, where $\Delta_1(X) = X^4 - X^2 +1$ and
$$\Delta_2(X) = 3X^4 - 2X^3 - X^2 - 2X + 3;$$ we have $\rho(\Delta) = 8$. The corresponding polynomial
$$P(X) = (-1)^4X^{8} \Delta(1-X^{-1})$$ is the one of example \ref{first example} : it is equal to $f_1f_2$, where 
$$f_1(X) = X^4 - 2X^3 + 5X^2 - 4X + 1,$$
and $$f_1(X) = X^4 - 2X^3 + 11X^2 - 10X + 3.$$ 

\medskip We have $\Pi_{f_1,f_2} = \{2\}$, hence $G_{\Delta} = G_P = 0$. 
Corollary \ref{knot coro signature} implies that there exists an $m$-knot with Alexander polynomial
$\Delta$ and signature $8$. But such a knot is indecomposable, since an $m$-knot with Alexander polynomial $\Delta_i$ has signature $0$ for $i = 1,2$.

\end{example} 

\begin{example}\label{6} Let $a \geqslant 0$ be an integer, and set $$\Delta_a(X) = X^6 -a X^5 - X^4 +(2a-1) X^3 - X^2 -a X + 1.$$ The polynomial $\Delta_a$ is irreducible, and $\rho(\Delta_a) = 4$ (see \cite {GM}, \S 7.3, Example 1 on page 284). This implies that all $m$-knots with Alexander polynomial $\Delta_a$ have signature 0. 

\medskip Let $b \geqslant 0$ be an integer with $b \not = a$. We have $\rho(\Delta_a \Delta_b) = 8$, and if 
moreover  $G_{\Delta_a \Delta_b} = 0$, then there exist $m$-knots with Alexander polynomial $\Delta_a \Delta_b$ and signature $8$;
these knots are indecomposable. We can take for instance $a = 0$ and $b = 2$; then $\Pi_{\Delta_a,\Delta_b} = \{2 \}$, hence $G_{\Delta_a \Delta_b} = 0$.

\end{example}

\section{$3$-knots in the $5$-sphere}\label{3}

The signature of a $3$-dimensional knot $K^3 \subset S^5$ is {\it divisible by $16$} (see for instance 
\cite{KW}, \S 3, page 95). 
The aim of this section is to show that with this additional restriction, the results of \S \ref{knot section}
extend to $3$-knots.

\medskip

Let $n \geqslant 1$ be an integer, and let $\Delta \in {\bf Z}[X]$ be a polynomial of degree $2n$ such
that $\Delta(X) = X^{2n}\Delta(X^{-1})$, $\Delta(1) = (-1)^n$ and that $\Delta(-1)$ is a square of an
integer. Set $P(X) = (-1)^nX^{2n} \Delta(1-X^{-1})$. Assume that $P$ is a product
of distinct irreducible monic polynomials $f \in {\bf Z}[X]$ such that $f(1-X) = f(X)$.
Let $G_{\Delta} = G_P$ be the group introduced in \S \ref{group}.

\medskip Let $s$ be an integer such that $s \equiv 0 \ {\rm (mod \ 16)}$, and that $|s| \leqslant \rho(P)$. 
Let $\tau \in {\rm Mil}_s(P)$. 

\begin{theo}\label{knot 3} There exists a $3$-knot with Alexander polynomial $\Delta$
and Milnor signature $\tau$ if and only if $\epsilon_{\tau} = 0$. 

\end{theo} 

\noindent
{\bf Proof.} This follows from Theorem \ref{main Seifert} and from a result of Levine (see \cite {Le 70}, Theorem 2) :
if $A$ is a Seifert form of signature divisible by $16$, then there exists a $3$-knot in the $5$-sphere
with Seifert form $S$-equivalent to $A$. Since $S$-equivalent Seifert forms have the same Alexander polynomial and Milnor signature, this completes the proof of the theorem. 

\begin{coro}\label{knot coro 3} Assume that $G_{\Delta} = 0$. Then for all $\tau \in {\rm Mil}_s(\Delta)$ there
exists a $3$-knot with Alexander polynomial $\Delta$ and Milnor signature $\tau$.

\end{coro}

\begin{coro}\label{knot coro signature}  Assume that $G_{\Delta} = 0$. Then there
exists a $3$-knot with Alexander polynomial $\Delta$ and signature $s$.

\end{coro}

\section{Unimodular Seifert forms}\label{unimodular}

We conclude by some remarks on a special case, which was already treated in detail in \cite {B 21}. Let $A : L \times L \to {\bf Z}$ be a unimodular Seifert form, i.e. ${\rm det}(A) = \pm 1$, and let $S : L \times L \to {\bf Z}$, defined by 
$S(x,y) =  A(x,y) + A(y,x)$, be the associated even, unimodular lattice. 

\medskip

Let
$t : L \to L$ be defined by $A(tx,y) = - A(y,x)$ for all $x,y \in L$; note that $t$ is an isometry
of $A$, and hence of $S$, and that the characteristic polynomial of $t$ is ${\rm det}(A) \Delta_A$.

\begin{prop} Sending $(L,A)$ to $(L,S,t)$ induces a bijection between isomorphism classes
of unimodular Seifert forms and isomorphism classes of even, unimodular lattices with an isometry.

\end{prop}

Hence the existence of a unimodular Seifert form with a given Alexander polynomial and Milnor signature
is equivalent to the existence of an even, unimodular lattice having an isometry of a given characteristic
polynomial and Milnor signature. This question is 
treated in \cite{B 21},  \S 25, 27 and 31; we recover the results of \S \ref{Seifert form section} in this
special case.

\bigskip
\bigskip
Eva Bayer--Fluckiger 

\medskip

EPFL-FSB-MATH

Station 8

1015 Lausanne, Switzerland

\medskip

eva.bayer@epfl.ch

\end{document}